\documentclass[a4paper,12pt]{article}
\usepackage{latexsym, amsmath, amsfonts, amscd, amssymb, verbatim, graphicx, graphics, fancyhdr}
\usepackage{color}
\usepackage{xcolor}

\def\tto{\;{\lower 1pt \hbox{$\rightarrow$}}\kern -10pt
\hbox{\raise 2pt \hbox{$\rightarrow$}}\;}

\def\proof{{\bf Proof.}}
\begin{document}
\pagestyle{myheadings}
\newtheorem{Theorem}{Theorem}[section]
\newtheorem{Proposition}[Theorem]{Proposition}
\newtheorem{Remark}[Theorem]{Remark}
\newtheorem{Lemma}[Theorem]{Lemma}
\newtheorem{Corollary}[Theorem]{Corollary}
\newtheorem{Definition}[Theorem]{Definition}
\newtheorem{Example}[Theorem]{Example}
\renewcommand{\theequation}{\thesection.\arabic{equation}}
\normalsize
\setcounter{equation}{0}

\title{\bf Approximate Proper Efficiency in Vector Optimization via Benson's Approach}

\author{Nguyen Thi Thu Huong\footnote{Center for Applied Mathematics and Informatics, Institute of Information and Communication Technology, Le Quy Don Technical University, 236 Hoang Quoc Viet Road, Bac Tu Liem District, Hanoi, Vietnam; email: nguyenhuong2308.mta@gmail.com.}}\maketitle

\date{}%\small\today                 B either coincides with M or is an empty set

\medskip
\begin{quote}
\noindent {\bf Abstract.} 
We present two criteria for checking approximate proper efficiency in vector optimization problems with the ordering cone being a nonnegative orthant. Although the criteria can be established by Benson's approach [H.P. Benson, \textit{An improved definition of proper efficiency for vector maximization with respect to cones}, J. Math. Anal. Appl. \textbf{71} (1979), 232--241], detailed proofs are given for the first time here. The two criteria are strong motivations to introduce the concept of $e$-properly efficient solution, where $e$ is any nonzero vector taken from the closed pointed convex ordering cone. For an arbitrary linear vector optimization problem, we show that either the $e$-properly efficient solution set is empty or it coincides with the $e$-efficient solution set. This new result has no analogue in the literature.

\noindent{\bf Mathematics Subject Classification (2010).}\ 90C26,
90C29, 90C33, 

\noindent {\bf Key Words.}\ Vector optimization; Linear vector optimization; pointed closed convex cone; proper efficiency; $\varepsilon$-proper efficiency; $e$-proper efficiency.
\end{quote}

\section{Introduction}
\markboth{Benson's Characterization for $\varepsilon$-proper Efficiency and Applications}{\sc n. t. t. huong}
\setcounter{equation}{0}

The concept of \textit{efficient solution} has played a fundamental role in analyzing vector optimization problems (see, e.g., \cite{Luc89,S86}). In general, the efficient solution set may include some points of a certain anomalous type.  By excluding unfavorable performances of the ratio between the profit rate of one cost and the corresponding loss of another, Geoffrion introduced in his pioneering work~\cite{GE68} the concept of \textit{properly  efficient solution}, where the ordering cone is the nonnegative orthant of the image space. Then,  Borwein~\cite{Bo77}, Benson~\cite{Be79}, Henig~\cite{Henig_JOTA1982} and other scholars have proposed significant extensions of Geoffrion’s concept of proper efficiency, which work for vector optimization problems where the ordering cone can be any pointed closed convex cone. The relationships among different types of properly efficient solutions can be found in the papers of Benson~\cite{Be79} and of Guerraggio et al.~\cite{GMZ94}. 

\medskip
When the ordering cone is the nonnegative orthant, Benson’s concept of proper efficiency~\cite{Be79}, which coincides with that of Geoffrion, deserves special attention. Isermann~\cite{I74} proved that Geoffrion's properly efficient solution set of a linear vector optimization problem coincides with the efficient solution set. Later, Choo~\cite{Choo84} extended this result to linear fractional vector optimization problems with bounded constraint sets. Recently, by using the recession cone of the constraint set and some properties of linear fractional functions (see \cite{LTY05,S86}), Huong et al.~\cite{HYY_JOGO2020,HYY2020} have obtained verifiable sufficient conditions for an efficient solution of a linear fractional vector optimization problem, where the constraint set is unbounded, to be a Geoffrion properly efficient solution. 

\medskip
\textit{Improperly efficient solutions} in the sense of Geoffrion in linear fractional vector optimization problems with unbounded constraint sets were studied systematically for the first time by Huong and Yen in~\cite{HY2022a}. Two sets of conditions assuring that all the efficient solutions of a given problem are improperly efficient have been given there. In addition, necessary conditions for an efficient solution to be improperly efficient and new sufficient conditions for Geoffrion's proper efficiency have been obtained in~\cite{HY2022a}.

\medskip
To satisfy some requirements of decision-makers one only needs to find \textit{approximate solutions} to the problems, or more often when there are no exact solutions to optimization problems. Therefore, considering approximate solutions and studying their necessary and sufficient conditions is an important issue from both theoretical and practical points of view. Different concepts of approximate solutions of vector optimization problems and the related optimality conditions can be found in~\cite{CK2016,GJN06,LiWang,ST14,Tammer94,White86} and the references therein. 

\medskip
Based on the idea of Geoffrion’s proper efficiency, Li and Wang~\cite{LiWang} introduced the concept of \textit{$\varepsilon$-proper efficiency} and obtained several necessary and sufficient conditions for $\varepsilon$-proper efficiency via scalarization and an alternative theorem. Then, for the approximate properly efficient solutions, Liu~\cite{Liu99} obtained some scalarization results. Very recently, Tuyen~\cite[Theorem~3.2]{T22} has shown that there is no difference between the $\varepsilon$-properly efficient solution set and the $\varepsilon$-efficient solution set for a linear fractional vector optimization problem with a bounded constraint set. But this fact does not hold even for linear vector optimization problems with unbounded constraint sets; see~\cite[Example~3.1]{T22}. 

\medskip
Recall that Benson's original characterizations~\cite[Theorem~3.2]{Be79} for the proper efficient solutions in the sense of Geoffrion allowed one to avoid using gain-to-loss ratios~\cite[p.~624]{GE68}, thus one could have an adequate extension of Geoffrion's notion of proper efficiency to vector optimization problems where the orderings are given by any nontrivial closed convex cones (see~\cite[Definition~2.4]{Be79}). By combining the approximate efficiency notion and Benson's approach \cite{Be79}, Rong and Ma \cite{RM2000} introduced the $e$-Benson proper efficiency concept where $e$ is taken from the ordering cone. Later, Guti\'errez et al. \cite{Gutierrez-et-al-12} and Zhao and Yang~\cite{ZhaoYang15} also proposed some concepts of approximate proper efficiency in the literature. Zhao et al. summarized several kinds of those concepts and discussed some relations between them in \cite{ZhaoChenYang15}. 

\medskip
The first aim of this paper is to establish two criteria for checking approximate proper efficiency in the sense of Li and Wang \cite{LiWang} as well as that of Liu \cite{Liu99} via Benson's approach~\cite{Be79}. To obtain the criteria, we use an extended version of Lemma~2.8 in~\cite{HWYY2019}, as well as some techniques from the classical papers Borwein~\cite{Bo77} and Benson~\cite{Be79}. Note that the first criterion was stated in~\cite[Theorem~4.1]{ZhaoChenYang15} without proof. The second criterion, where the pointedness of the cone $\mathbb R^m_+$ is taken into account, is a refined version of the first one. The second aim of this paper is to prove that, for any linear vector optimization problem with a pointed polyhedral convex cone $K$ and for any $e\in K\setminus\{0\}$, either the $e$-properly efficient solution set is empty or it coincides with the $e$-efficient solution set. For $e=0$, it is well-known~\cite{I74} that the $e$-properly efficient solution set of such a problem coincides with the $e$-efficient solution set.

\medskip
The paper is organized as follows. Section~\ref{Sect_2} recalls some notions, definitions, and auxiliary results. Section~\ref{Sect_3} presents the neccessary and sufficient conditions for $\varepsilon$-proper efficiency of a general vector optimization problem. Section~\ref{Sect_4} gives several illustrative examples. Some interesting results on $e$-proper efficiency in linear vector optimization are established in Sect.~\ref{Sect_5}.

\section{Preliminaries}\label{Sect_2}
\markboth{Approximate Proper Efficiency in Vector Optimization}{\sc n. t. t. huong}
\setcounter{equation}{0}
The scalar product and the norm in the Euclidean space $\mathbb R^n$ are denoted, respectively, by $\langle\cdot,\cdot\rangle$ and $\|\cdot\|$. 
Vectors in $\mathbb R^n$ are represented by columns of real numbers. If $A$ is a matrix, then $A^T$ denotes the transposed matrix of $A$. Thus, one has $\langle x,y\rangle=x^Ty$  for any $x,y\in \mathbb R^n$. The nonnegative orthant in $\mathbb R^n$ and the set of positive integers are denoted respectively by $\mathbb R^n_+$ and $\mathbb N$. 

A nonempty set $K\subset\mathbb R^m$ is called a \textit{cone} if $tv\in K$ for all $v\in K$ and $t\geq 0$. One says that $K$ is \textit{pointed} if $K\cap (-K)=\{0\}$. The smallest cone containing a nonempty set $D\subset\mathbb R^m$, i.e., the cone generated by $D$, will be denoted by ${\rm cone}\,D$. The topological closure of $D$ is denoted by $\overline{D}$ and $\overline{\rm cone}\,D:=\overline{{\rm cone}\, D}$.

A nonzero vector $v\in\mathbb R^n$ (see \cite[p.~61]{Roc70}) is said to be a \textit{direction of recession} of a nonempty convex set $D\subset\mathbb R^n$ if $x+tv\in D$ for every $t\geq 0$ and every $x\in D.$  The set composed by $0\in\mathbb R^n$ and all the directions $v\in\mathbb R^n\setminus\{0\}$ satisfying the last condition, is called the \textit{recession cone} of $D$ and denoted by $0^+D.$  If $D$ is closed and convex, then  $$0^+D=\{v\in\mathbb R^n\;:\; \exists x\in\Omega\ \, {\rm s.t.}\ \, x+tv\in D\ \, {\rm for\ all}\ \, t> 0\}.$$

Let $K\subset\mathbb R^m$ be a nonempty closed convex cone. Let there be given a nonempty subset $X\subset \mathbb R^n$ and functions $f_i:X\to\mathbb R$ with $i\in I$, where $I:=\{1,\cdots,m\}$. Put $$f(x)=(f_1(x),\dots,f_m(x))^T$$ for $x\in X.$ For any $y^2,y^1\in\mathbb R^m$, if $y^2-y^1\in K$, one writes $y^1\leq_K y^2$ and also $y^2\geq_K y^1$. Consider the vector optimization problem 

\medskip
\hskip1cm (VP) \hskip 2cm ${\rm Minimize}\ \, f(x)\ \,$
subject\ to $\; x\in X$.

\medskip
If $x\in X$ and there exists no $y\in X$ such that
$f(y)\leq_K f(x)$ and $f(y)\neq f(x)$, then one says that $x$ is an {\it efficient solution} of (VP). The efficient solution set of~(VP) is denoted by $E$. 

One says that (VP) is a \textit{linear vector optimization problem} if all the functions $f_i$ for $i\in I$ are  affine, $X$ is polyhedral convex set (see, e.g.,~\cite[Section~19]{Roc70}), and $K$ is a polyhedral convex cone.

Let $\varepsilon=(\varepsilon_1,...,\varepsilon_m)$ be a vector in $\mathbb{R}^m_+$. 

\begin{Definition}\label{appro_def_efficiency}{\rm (See~\cite{LiWang,Liu99}) A point $x\in X$ is said to be an {\it $\varepsilon$-efficient solution} of (VP) if there exists no $y\in X$ such that $f(y)\leq_K f(x)-\varepsilon$ and $f(y)\neq f(x)-\varepsilon$.} 
\end{Definition}

The $\varepsilon$-efficient solution set of (VP) is denoted by $E_{\varepsilon}$.

\begin{Remark}\label{Remark} {\rm Let $\bar x \in X$. Then, $\bar x$ is an $\varepsilon$-efficient solution of (VP) if and only if ${\rm cone}[f(X)-(f(\bar x)-\varepsilon)]\cap (-K)=\{0\}$. If $K$ is pointed, it is easy to show that the latter is equivalent to the condition  ${\rm cone}[f(X)+K-(f(\bar x)-\varepsilon)]\cap (-K)=\{0\}$.}
\end{Remark}

When $\varepsilon=0$, the notion of $\varepsilon$-efficient solution reduces to the notion of efficient solution, i.e., $E_0=E$.

Geoffrion’s definition of a properly efficient solution~\cite[p.~618]{GE68} applies to the case $K=\mathbb R^m_+$, where $\mathbb R^m_+$ denotes the nonnegative orthant of $\mathbb R^m$. Namely, $\bar x\in E$ is said to be a \textit{Geoffrion's properly efficient solution} of~(VP) if there exists a positive constant $M$ such that, for each $i\in I$, whenever $x\in K$ and $f_i(x)<f_i(\bar x)$ one can find an index $j\in I$ such that $f_j(x)>f_j(\bar x)$ and $A_{i,j}(\bar x,x)\leq M$ with $A_{i,j}(\bar x,x):=\displaystyle\frac{f_i(\bar x)-f_i(x)}{f_j(x)-f_j(\bar x)}$. Geoffrion's efficient solution set of (VP) is denoted by $E^{Ge}$.

\begin{Definition}\label{Def appro Proper efficient solution} {\rm (See~\cite{LiWang,Liu99}) One says that $\bar x\in E_\varepsilon$ is an \textit{$\varepsilon$-properly efficient solution} of (VP), where $K=\mathbb R^m_+$, if there exists a scalar $M>0$ such that, for each $i\in I$, whenever $x\in X$ and $f_i(x)<f_i(\bar x)-\varepsilon_i$ one can find an index $j\in I$ such that $f_j(x)>f_j(\bar x)-\varepsilon_j$ and $A_{i,j}(\bar x,x,\varepsilon)\leq M$ with $A_{i,j}(\bar x,x,\varepsilon):=\displaystyle\frac{f_i(\bar x)-f_i(x)-\varepsilon_i}{f_j(x)-f_j(\bar x)+\varepsilon_j}$.}
\end{Definition}

In the case $K=\mathbb R^m_+$, the set of the above-defined $\varepsilon$-properly efficient solutions of (VP) is abbreviated to $E^{Ge}_{\varepsilon}$. Note that $E^{Ge}_0=E^{Ge}$.

\begin{Remark}\label{Rem1}
	{\rm (See~\cite[p.~110]{Liu99}) For any $\bar x\in E_{\varepsilon}$, $\bar x\notin E^{Ge}_{\varepsilon}$ if and only if for every scalar $M>0$ there exist $x\in X$ and $i\in I$ with $f_i(x)<f_i(\bar x)-\varepsilon_i$ such that, for all $j\in I$ satisfying $f_j(x)>f_j(\bar x)-\varepsilon_j$, one has $A_{i,j}(\bar x,x,\varepsilon)>M$.}
\end{Remark}

The following lemma is an extension of Lemma~2.8 in~\cite{HWYY2019}, where the case $K=\mathbb R^m_+$ was considered. 

\begin{Lemma}\label{equivalence} Let $K\subset\mathbb R^m$ be a pointed closed convex cone. Then, for any nonempty subset $D\subset\mathbb R^m$, one has
	\begin{equation}\label{intersection1}
		\overline{\rm cone}\,{\big(D+K\big)}\cap \big({-K}\big)=\{0\}
	\end{equation}
	if and only if 
	\begin{equation}\label{intersection2}
		\overline{\rm cone}\,D\cap \big({-K}\big)=\{0\}.
	\end{equation}
\end{Lemma}
\proof\ Since $0\in K$ by the pointedness of the cone, one has $D\subset D+K$. So, the ``only if" assertion is clear. To prove the ``if" part of the lemma, suppose on the contrary that~\eqref{intersection2} holds, but~\eqref{intersection1} is invalid. Then, there are a nonzero vector $v\in-K$, a sequence $\{t_k\}$ of positive real numbers, and sequences $\{r^k\}\subset K$, $\{a^k\}\subset D$, such that $v=\displaystyle\lim_{k\to\infty}\big[t_k(a^k+r^k)\big].$ Setting $u^k=t_kr^k$ for all $k\in\mathbb N$, we have $\{u^k\}\subset K$ and \begin{equation}\label{expression_for_v} v=\displaystyle\lim_{k\to\infty}\big(t_ka^k+u^k\big).\end{equation}

If the sequence $\{u^k\}$ is bounded, we may assume that  $\displaystyle\lim_{k\to\infty}u^k=u$. By the closedness of $K$, $u\in K$. From~\eqref{expression_for_v} it follows that $$\displaystyle\lim_{k\to\infty}\left(t_ka^k\right)=v-u\in- K\setminus\{0\}-K.$$
Hence, by the pointedness and the convexity of the cone $K$, we have $$\displaystyle\lim_{k\to\infty}\left(t_ka^k\right)\in -K\setminus\{0\}.$$ Since $\displaystyle\lim_{k\to\infty}\left(t_ka^k\right)\in \overline{\rm cone}\,D$, this contradicts~\eqref{intersection2}.

If $\{u^k\}$ is unbounded, we may assume that $\displaystyle\lim_{k\rightarrow\infty}\|u^k\|=+\infty$, $u^k\neq 0$ for all~$k$, and $\displaystyle\lim_{k\rightarrow\infty}\frac{u^k}{\|u^k\|}=z$, where $\|z\|=1$. By the closedness of $K$, $z\in K$. Hence $-z\in -K$. From~\eqref{expression_for_v} it follows that
$$\displaystyle\lim_{k\rightarrow\infty}\dfrac{v}{\|u^k\|}=\displaystyle\lim_{k\rightarrow\infty}\Big(\dfrac{t_k}{\|u^k\|}a^k+\dfrac{u^k}{\|u^k\|}\Big).$$ 
Passing the last expression to limit as $k\rightarrow\infty$ and note that $\displaystyle\lim_{k\rightarrow\infty}\frac{u^k}{\|u^k\|}=z$, one gets
$$-z=\displaystyle\lim_{k\rightarrow\infty}\Big(\dfrac{t_k}{\|u^k\|}a^k\Big)\in\overline{\rm cone}\,D.$$ Thus $0\ne -z\in\overline{\rm cone}\,D\cap \big({-K}\big)$. This contradicts \eqref{intersection2} and completes the proof. $\hfill\Box$

\section{Criteria for $\varepsilon$-Proper Efficiency}~\label{Sect_3}
\markboth{Approximate Proper Efficiency in Vector Optimization}{\sc n. t. t. huong}
\setcounter{equation}{0}

To obtain a criterion for $\varepsilon$-proper efficiency, we will establish two propositions. The proof of the first proposition relies on some ideas of Borwein~\cite{Bo77}. Meanwhile, the proof of the second one uses several arguments of Benson~\cite{Be79}.  

\begin{Proposition}\label{necessary condition} {\rm (A neccessary condition for $\varepsilon$-proper efficiency)} If $\bar x\in X$ is an $\varepsilon$-properly efficient solution of {\rm (VP)} in the case $K=\mathbb R^m_+$, then 	
	\begin{equation}\label{appro_Be-def}
		\overline{\rm cone}\big[f(X)+\mathbb R^m_+-(f(\bar x)-\varepsilon)\big]\cap \big({-\mathbb R^m_+}\big)=\{0\}.
	\end{equation}
\end{Proposition}
\proof\ Let $\bar x\in E^{Ge}_{\varepsilon}$. Clearly,
$$0\in \overline{\rm cone}\big[f(X)+\mathbb R^m_+-(f(\bar x)-\varepsilon)\big]\cap \big({-\mathbb R^m_+}\big).$$ So, if~\eqref{appro_Be-def} fails to hold, then there is a vector $v=(v_1,...,v_m)$ belonging to the set on the left-hand side of~\eqref{appro_Be-def} with $v\ne 0$. 
Hence, we can find sequences $\{x^k\}\subset X$, $\{u^k\}\subset\mathbb R^m_+$, and $\{\tau_k\}\subset\mathbb R_+$  such that $\displaystyle\lim_{k\to\infty}\tau_k\big(f(x^k)+u^k-f(\bar x)+\varepsilon\big)=v$. Since $v\in\big({-\mathbb R^m_+}\big)\setminus\{0\}$, we see that  $v_i\leq 0$ for all $i\in I$ and there exists $i_0\in I$ such that $v_{i_0}<0$. As $\displaystyle\lim_{k\to\infty}\tau_k\big(f_{i_0}(x^k)+u^k_{i_0}-f_{i_0}(\bar x)+\varepsilon_{i_0}\big)=v_{i_0}<0$, there is $k_0\in\mathbb N$ with  \begin{equation}\label{ineq_k}\tau_k\big(f_{i_0}(x^k)+u^k_{i_0}-f_{i_0}(\bar x)+\varepsilon_{i_0}\big)<\dfrac{v_{i_0}} {2}
\end{equation} for all $k\geq k_0$. By~\eqref{ineq_k} one has \begin{equation}\label{ineq_k_new} f_{i_0}(x^k)<f_{i_0}(\bar x)-\varepsilon_{i_0}\quad (\forall k\geq k_0).\end{equation} 
Since $\bar x$ is an $\varepsilon$-efficient solution of~(VP), this implies that the set $$I_k:=\{j\in I: f_j(x^k)>f_j(\bar x)-\varepsilon_j\}$$ is nonempty for every $k\geq k_0$. Therefore, by working with a subsequence if necessary, we may assume that $I_k=\widetilde I$ for all $k\geq k_0$, where $\widetilde I$ is a nonempty subset of $I$. Hence, for every $k>k_0$, by~\eqref{ineq_k} we have $\tau_k>0$ and
\begin{equation}\label{eq i_0}
	f_{i_0}(\bar x)-f_{i_0}(x^k)-\varepsilon_{i_0}>\dfrac{|v_{i_0}|}{2\tau_k}.
\end{equation} Let $M>0$ be given arbitrarily. Since  $\displaystyle\lim_{k\to\infty}\tau_k\big(f_{i}(x^k)-f_{i}(\bar x)+\varepsilon_{i}+u^k_{i}\big)\leq 0$ for all $i\in I$, then there exists an integer $k_M$ with $k_M\geq k_0$ such that 
\begin{equation}\label{eq i}
	\tau_k\big(f_{i}(x^k)-f_{i}(\bar x)+\varepsilon_{i}+u^k_{i}\big)<\dfrac{|v_{i_0}|}{2M}  
\end{equation}
for any $k\geq k_M$ and $i\in I$. By~\eqref{eq i},
$$f_i(x^k)-f_i(\bar x)+\varepsilon_i <\dfrac{|v_{i_0}|}{2\tau_kM}$$ for any $k\geq k_M$ and $i\in I$.
Thus, for every $j\in\widetilde I$ one has 
\begin{equation}\label{ineq j}
	0<f_j(x^k)-f_j(\bar x)+\varepsilon_j <\dfrac{|v_{i_0}|}{2\tau_kM}\quad (\forall k>k_M).
\end{equation}
Therefore, from \eqref{eq i_0} and \eqref{ineq j}, it follows that
$$A_{i_0,j}(\bar x,x^k,\varepsilon)=\displaystyle\frac{f_{i_0}(\bar x)-f_{i_0}(x^k)-\varepsilon_{i_0}}{f_j(x^k)-f_j(\bar x)+\varepsilon_j}>\displaystyle\frac{\dfrac{|v_{i_0}|}{2\tau_k}}{\dfrac{|v_{i_0}|}{2\tau_kM}}=M$$ for all $k\geq k_M$ and for every $j\in\widetilde I=I_k$. Thus, for any $M>0$ there exist $x^k\in X$ and  $i_0\in I$ with $f_{i_0}(x)<f_{i_0}(\bar x)-\varepsilon_{i_0}$ (see~\eqref{ineq_k_new}) satisfying $A_{i_0,j}(\bar x,x^k,\varepsilon)>M$ for all $j\in I$ with $f_j(x^k)>f_j(\bar x)-\varepsilon_j$. Hence, by Remark~\ref{Rem1}, $\bar x\notin E^{Ge}_{\varepsilon}$. We have arrived at a contradiction, which completes the proof.
$\hfill\Box$

\begin{Proposition}\label{sufi condiction} {\rm (A sufficient condition for $\varepsilon$-proper efficiency)}
	A point $\bar x\in X$ is an $\varepsilon$-properly efficient solution of {\rm (VP)} in the case $K=\mathbb R^m_+$ if the equality~\eqref{appro_Be-def} holds.	                                          
\end{Proposition} 
\proof\ Suppose that~\eqref{appro_Be-def} holds. Clearly, from~\eqref{appro_Be-def} it follows that 
$$\big[f(X)-(f(\bar x)-\varepsilon)\big]\cap \big({-\mathbb R^m_+}\big)\subset\{0\}.$$ Hence, in accordance with Definition~\ref{appro_def_efficiency}, one has $\bar x\in E_{\varepsilon}$.

Suppose on the contrary that~$\bar x\notin E^{Ge}_{\varepsilon}$. Then, by Remark~\ref{Rem1}, for every $k\in\mathbb N$ we can find  $x^k\in K$ and $i(k)\in I$ with $f_{i(k)}(x^k)<f_{i(k)}(\bar x)-\varepsilon_i$ such that $A_{i(k),j}(\bar x,x^k,\varepsilon)>k$ for all $j\in I$ satisfying $f_j(x^k)>f_j(\bar x)-\varepsilon_j$. Since the sequence $\{i(k)\}$ has values in the finite set~$I$, by considering a subsequence, we may assume that $i(k)=i$ for all $k$ with $i\in I$ being a fixed index. For each $k$, as $\bar x\in E_\varepsilon$ and $f_{i}(x^k)<f_{i}(\bar x)-\varepsilon_i$, there exist an index $j(k)\in I\setminus\{i\}$ satisfying  $f_{j(k)}(x^k)>f_{j(k)}(\bar x)-\varepsilon_j$. Again, working with a subsequence if necessary, we may assume that $j(k)=j$ for all $k$, where $j\in I\setminus\{i\}$ is a fixed index. Thus,
\begin{equation}\label{basic_inequality}
	A_{i,j}(\bar x, x^k,\varepsilon)=\frac{f_i(\bar x)-f_i(x^k)-\varepsilon_i}{f_j(x^k)-f_j(\bar x)+\varepsilon_j}>k\quad\ (\forall k\in\mathbb N).
\end{equation}
In addition, we can assume that the index set $$\widetilde I:=\{j'\in I\setminus\{i\}\,|\,f_{j'}(x^k)>f_{j'}(\bar x)-\varepsilon_{j'}\}$$ is constant for all $k$. Note that $j\in\widetilde I$. For each $k$, putting 
\begin{equation}\label{lambda k}
	\lambda_k=\frac{1}{f_i(\bar x)-f_i(x^k)-\varepsilon_i},
\end{equation} one has $\lambda_k>0$. 

First, consider the case where $I\setminus\widetilde I=\{i\}$. For $\ell=i$ one has  
\begin{equation}\label{ell 1}\lambda_k\big(f_\ell(x^k)-f_\ell(\bar x)+\varepsilon_\ell\big) =-1\quad (\forall k\in\mathbb N).\end{equation}
For every $\ell\in\widetilde I$, from~\eqref{basic_inequality} we get 
\begin{equation}\label{inequation}
	0<\dfrac{1}{A_{i,\ell}(\bar x, x^k,\varepsilon)}<\dfrac{1}{k}\quad\ (\forall k\in\mathbb N).
\end{equation}
Passing the inequalities~\eqref{inequation} to limit as $k\rightarrow\infty$ and using~\eqref{lambda k} yield \begin{equation}\label{ell 2}\displaystyle\lim_{k\to\infty}\left[\lambda_k\big(f_\ell(x^k)-f_\ell(\bar x)+\varepsilon_\ell\big)\right] =\displaystyle\lim_{k\to\infty}\dfrac{1}{A_{i,\ell}(\bar x, x^k,\varepsilon)}=0.\end{equation} 
By~\eqref{ell 1} and~\eqref{ell 2}, the numbers 
\begin{equation}\label{h} h_\ell:=\displaystyle\lim_{k\to\infty}\left[\lambda_k\big(f_\ell(x^k)-f_\ell(\bar x)+\varepsilon_\ell\big)\right]\quad (\ell=1,\ldots,m)\end{equation} 
are well defined. Set $h=(h_1,\ldots,h_m)$ and observe that $h\in-\mathbb{R}^m_+\setminus\{0\}$. Since  $$\lambda_k\big(f(x^k)-f(\bar x)+\varepsilon\big)\in {\rm cone}\big[f(X)-(f(\bar x)-\varepsilon)\big]\quad\ (\forall k\in\mathbb N),$$ from~\eqref{h} it follows that $h\in\overline{\rm cone}\big[f(X)-(f(\bar x)-\varepsilon)\big]$. Thus, the nonzero vector $ h$ is contained in the set $\overline{\rm cone}\big[f(X)+\mathbb R^m_+-(f(\bar x)-\varepsilon)\big]\cap \big({-\mathbb R^m_+}\big).$ This contradicts~\eqref{appro_Be-def}.

Now, consider the case where $I\setminus\widetilde I\ne\{i\}$. As in the previous case,~\eqref{ell 1} is valid for $\ell=i$ and~\eqref{ell 2} holds for every $\ell\in\widetilde I$. For each $\ell\in I\setminus\widetilde I$ with $\ell\ne i$ and for every $ k\in\mathbb{N}$, we put  
\begin{equation}\label{v-k-ell} v^k_\ell=-\left[f_\ell(x^k)-f_\ell(\bar x)+\varepsilon_\ell\right].\end{equation} 
Since $\ell\notin \widetilde I$, one has $f_\ell(x^k)-f_\ell(\bar x)+\varepsilon_\ell\leq 0$; hence  $v^k_\ell\geq 0$. Next, choose $v^k_\ell=0$ for each $\ell\in\widetilde I\cup \{i\}$ and for every $k\in\mathbb{N}$. Then, for each $k\in\mathbb{N}$, the vector $v^k:=(v^k_1,\ldots,v^k_m)$ belongs to $\mathbb{R}^m_+$. Therefore,  \begin{equation}\label{inclusion_cone}\lambda_k\big(f(x^k)+v^k-f(\bar x)+\varepsilon\big)\in {\rm cone}\big[f(X)+\mathbb R^m_+-(f(\bar x)-\varepsilon)\big]\quad\ (\forall k\in\mathbb N).\end{equation} 
Define \begin{equation}\label{h_Caseb} h_\ell=\displaystyle\lim_{k\to\infty}\left[\lambda_k\big(f_\ell(x^k)+v^k_\ell-f_\ell(\bar x)+\varepsilon_\ell\big)\right]\quad (\ell=1,\ldots,m).\end{equation}  
Clearly,~\eqref{ell 1},~\eqref{ell 2},~\eqref{v-k-ell} and~\eqref{h_Caseb} imply that $h_i=-1$ and  $h_\ell=0$ for every $\ell\in I\setminus\{i\}$. In particular, $h\in -\mathbb R^m_+\setminus\{0\}$. So, by~\eqref{inclusion_cone} and~\eqref{h_Caseb} we can infer that 
$$\overline{\rm cone}\big[f(X)+\mathbb R^m_+-(f(\bar x)-\varepsilon)\big]\cap \big({-\mathbb R^m_+}\big)\neq \{0\}.$$ We have arrived at a contradiction.

The proof is complete. $\hfill\Box$

\begin{Theorem}\label{thm_criterion} {\rm (A criterion for $\varepsilon$-proper efficiency)}
	A point $\bar x\in X$ is an $\varepsilon$-properly efficient solution of {\rm (VP)} in the case $K=\mathbb R^m_+$ if and only if either~\eqref{appro_Be-def} or the equality
	\begin{equation}\label{appro_Be-def-s}
		\overline{\rm cone}\big[f(X)-(f(\bar x)-\varepsilon)\big]\cap \big({-\mathbb R^m_+}\big)=\{0\}
	\end{equation} holds.             	                        
\end{Theorem}
\proof\ Let $\bar x\in X$ be given arbitrarily. Then, by Propositions~\ref{necessary condition} and~\ref{sufi condiction}, $\bar x$ is an $\varepsilon$-properly efficient solution of {\rm (VP)} if and only if the equality~\eqref{appro_Be-def} is valid. Applying Lemma~\ref{equivalence} for $K=\mathbb R^m_+$ and $D=f(X)-(f(\bar x)-\varepsilon)$, we can infer that~\eqref{appro_Be-def} is valid if and only if~\eqref{appro_Be-def-s} holds. This justifies the assertion of the theorem. $\hfill\Box$

\medskip
Based on Definition~\ref{appro_def_efficiency}, Definition~\ref{Def appro Proper efficient solution}, and Theorem~\ref{thm_criterion}, we can extend the notions $\varepsilon$-efficient solution of $\varepsilon$-properly efficient solution in to vector optimization problems where the orderings are given by any closed convex cones as follows.

\medskip
Let $e$ be a vector in $K$. 

\begin{Definition}\label{appro_def_efficiency_e} {\rm A point $x\in X$ is said to be an {\it $e$-efficient solution} of (VP) if there exists no $y\in X$ such that $f(y)\leq_K f(x)-e$ and $f(y)\neq f(x)-e$.} 
\end{Definition}

The $e$-efficient solution set of (VP) is denoted by $E_e$. 

\begin{Definition}\label{Benson'appro pro_e} {\rm A point $\bar x\in X$ is said to be an \it{$e$-proper efficient solution}, where $e\in K\setminus \{0\}$, of (VP) if
		\begin{equation}\label{appro_Be-def_K} 
			\overline{\rm cone}\big[f(X)+K-(f(\bar x)-e)\big]\cap \big({-K}\big)=\{0\}.
	\end{equation}}
\end{Definition} 

The set of all $e$-properly efficient solutions of (VP) is denoted by $E^{Be}_e$. When $e=0$, the notion of $e$-properly efficient solution reduces to the notion of properly efficient solution which defined by Benson~\cite[Definition~2.4]{Be79}, i.e., $E^{Be}=E^{Be}_0$. Since~\eqref{appro_Be-def_K} surely yields $\big(f(X)-(f(\bar x)-e)\big)\cap \big({-K}\big)=\{0\}$, property \eqref{appro_Be-def_K} implies that $\bar x\in E_e$.

\begin{Remark}\label{equivalence e} 
	{\rm If $K$ is pointed, then one can rewrite~\eqref{appro_Be-def_K} equivalently as \begin{equation}\label{def_K_pointed} 
			\overline{\rm cone}\big[f(X)-(f(\bar x)-e)\big]\cap \big({-K}\big)=\{0\}.
		\end{equation} Indeed, applying Lemma~\ref{equivalence} for the nonemtyset set $D:=f(X)-(f(\bar x)-e).$}
\end{Remark}

Some results on the set $E^{Be}_e$ in linear vector optimization will be obtained in Section~\ref{Sect_5}. 

\section{Illustrative Examples}~\label{Sect_4}
\markboth{Approximate Proper Efficiency in Vector Optimization}{\sc n. t. t. huong}
\setcounter{equation}{0}

Our first example is an illustration for Theorem~\ref{thm_criterion}.

\begin{Example}{\rm (See~\cite{T22} and Figure~\ref{Fig.1})}\label{T_Exam1} {\rm Consider problem $({\rm VP})$ where $n=m=2$, $X=K=\mathbb{R}^2_+$, and $f(x)=(-x_1, x_2)$. Clearly, the problem has no efficient solution. Let $\varepsilon=(\varepsilon_1, \varepsilon_2)\in\mathbb{R}^2_+$ be such that $\varepsilon_2>0$. Using Remark~\ref {Remark},  it is easy to check that  
		$$E_{\varepsilon}=\{x=(x_1,x_2):x_1\geq 0, \, 0\leq x_2<\varepsilon_2\}.$$
		For each $\bar x\in E_{\varepsilon}$, we have $$\overline{\rm cone}\big[f(X)-(f(\bar x)-\varepsilon)\big]\cap \big({-\mathbb R^2_+}\big)=\{x=(x_1,0)\,:\,x_1\leq 0\}.$$ Thus, by Theorem~\ref{thm_criterion}, every point $\bar x\in E_{\varepsilon}$ is not an $\varepsilon$-properly efficient solution. So, $E^{Be}_\varepsilon=\emptyset$, while $E_{\varepsilon}\neq\emptyset$.}
	\begin{figure}[h]
		\begin{center}
			\includegraphics[height=5cm,width=8cm]{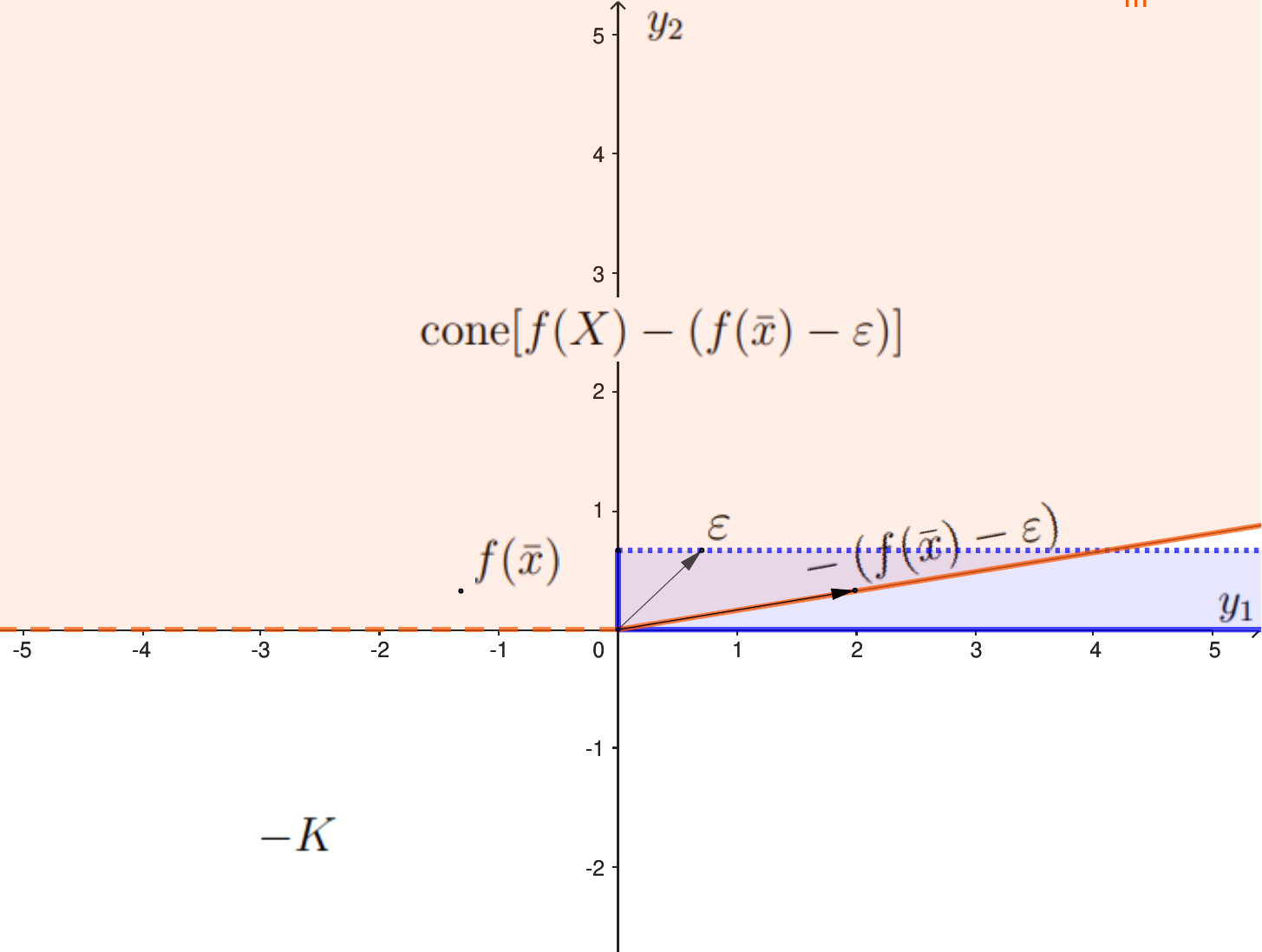}
			\caption{Illustration of Example \ref{T_Exam1}.}\label{Fig.1}
		\end{center}
	\end{figure}
\end{Example}	

The following example can be regarded as a simplified version of Example~\ref{T_Exam1}.

\begin{Example}\label{T_Exam2} {\rm (See Figure~\ref{Fig 2})} {\rm Consider problem $({\rm VP})$ with $n=m=2$, $K=\mathbb{R}^2_+$, $$X=\{x=(x_1,x_2)\in \mathbb{R}^2\,:\,x_1\geq 0\},$$ and $f(x)=(x_1, x_2)$. One has $E=\emptyset$. Let $\varepsilon=(\varepsilon_1, \varepsilon_2)\in\mathbb{R}^2_+$ be such that $\varepsilon_1>0$. Using Remark~\ref {Remark}, we can check that  $$E_{\varepsilon}=\{x=(x_1,x_2)\in\mathbb R^2\,:\, 0\leq x_1<\varepsilon_1\}.$$ Here, for every $\bar x\in E_{\varepsilon}$, we have $$\overline{\rm cone}\big[f(X)-(f(\bar x)-\varepsilon)\big]=\{x=(x_1,x_2)\,:\, x_1\geq 0,\,x_2\in\mathbb{R}\}.$$ So \eqref{appro_Be-def-s} does not hold. Then, by Theorem~\ref{thm_criterion}, every point $\bar x\in E_{\varepsilon}$ is not an $\varepsilon$-properly efficient solution. This means that $E^{Be}_\varepsilon=\emptyset$, while $E_\varepsilon\neq\emptyset$.
		\begin{figure}[h]
			\begin{center}
				\includegraphics[height=6cm,width=9cm]{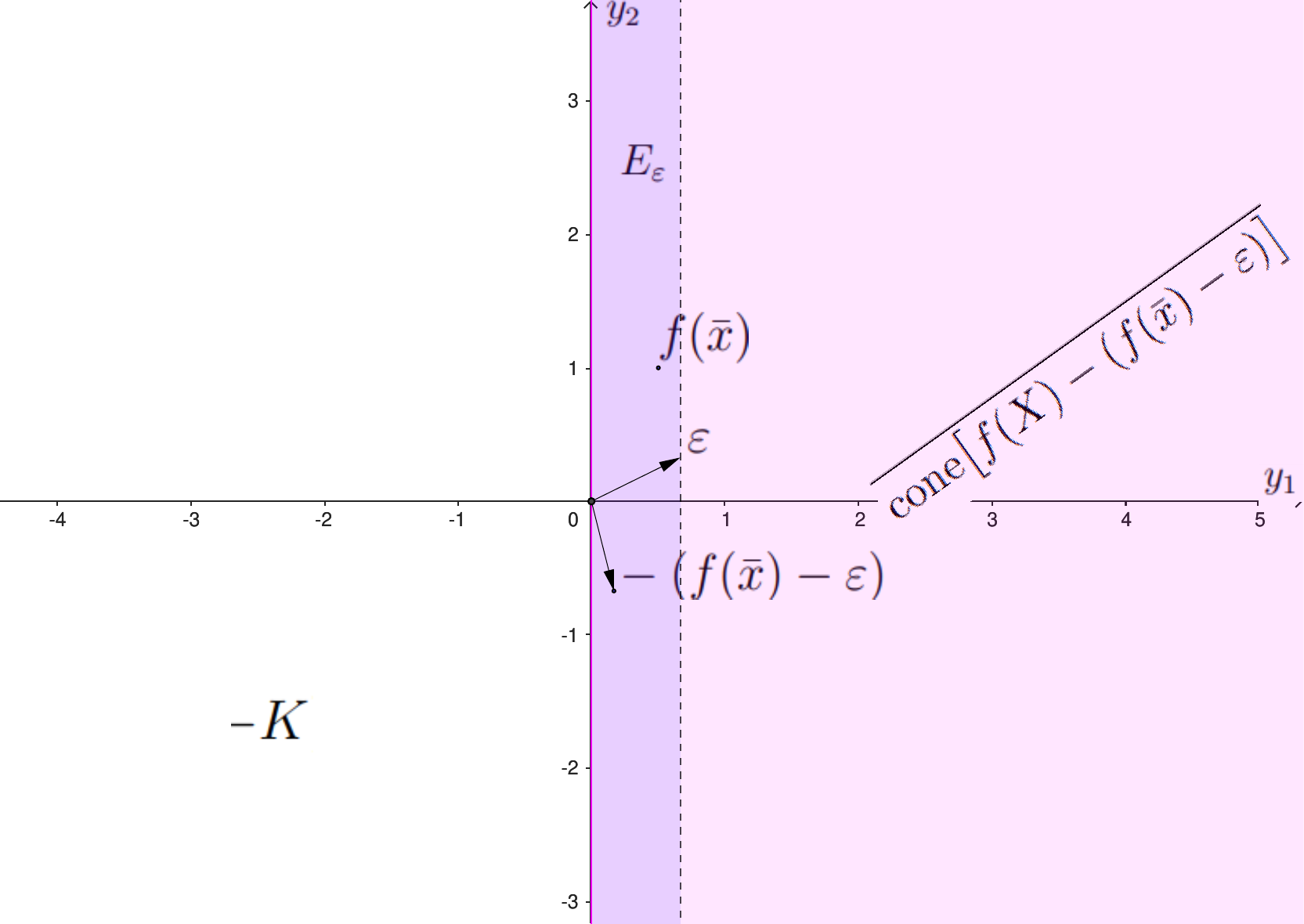}
				\caption{Illustration of Example \ref{T_Exam2}.}\label{Fig 2}
			\end{center}
	\end{figure}}  
\end{Example}

In the next example, one has $E_\varepsilon=E^{Be}_\varepsilon$ for every $\varepsilon\in\mathbb{R}^m_+$.

\begin{Example}\label{T_Exam3}{\rm (See Figure~\ref{Fig. 3})} {\rm Consider problem $({\rm VP})$ where $n=m=2$, $$K=\{v=(v_1,v_2)\in \mathbb{R}^2\,:\,0\leq v_2\leq v_1\},$$ $X=\{x=(x_1,x_2)\in \mathbb{R}^2\,:\,x_1\geq 0\}$, and $f(x)=(x_1, x_2)$. We can see that $$E=\{x=(x_1,x_2)\in\mathbb{R}^2\,:\, x_1=0 \}.$$ Let $e=(e_1, e_2)\in K\subset\mathbb{R}^2_+$. Using Remark~\ref {Remark}, we can check that  $$E_e=\{x=(x_1,x_2)\in\mathbb R^2\,:\, 0\leq x_1\leq e_1\}.$$ Here, for every $\bar x\in E_e$, we have $$\overline{\rm cone}\big[f(X)-(f(\bar x)-e)\big]=\{x=(x_1,x_2)\,:\, x_1\geq 0,\,x_2\in\mathbb{R}\}.$$ Thus, \eqref {def_K_pointed} is true. Then, by Definition~\ref{Benson'appro pro_e}, we have $E_e=E^{Be}_e$. In particular, taking $e=0$ gives $E^{Be}=E$.
		\begin{figure}[h]
			\begin{center}
				\includegraphics[height=6cm,width=9cm]{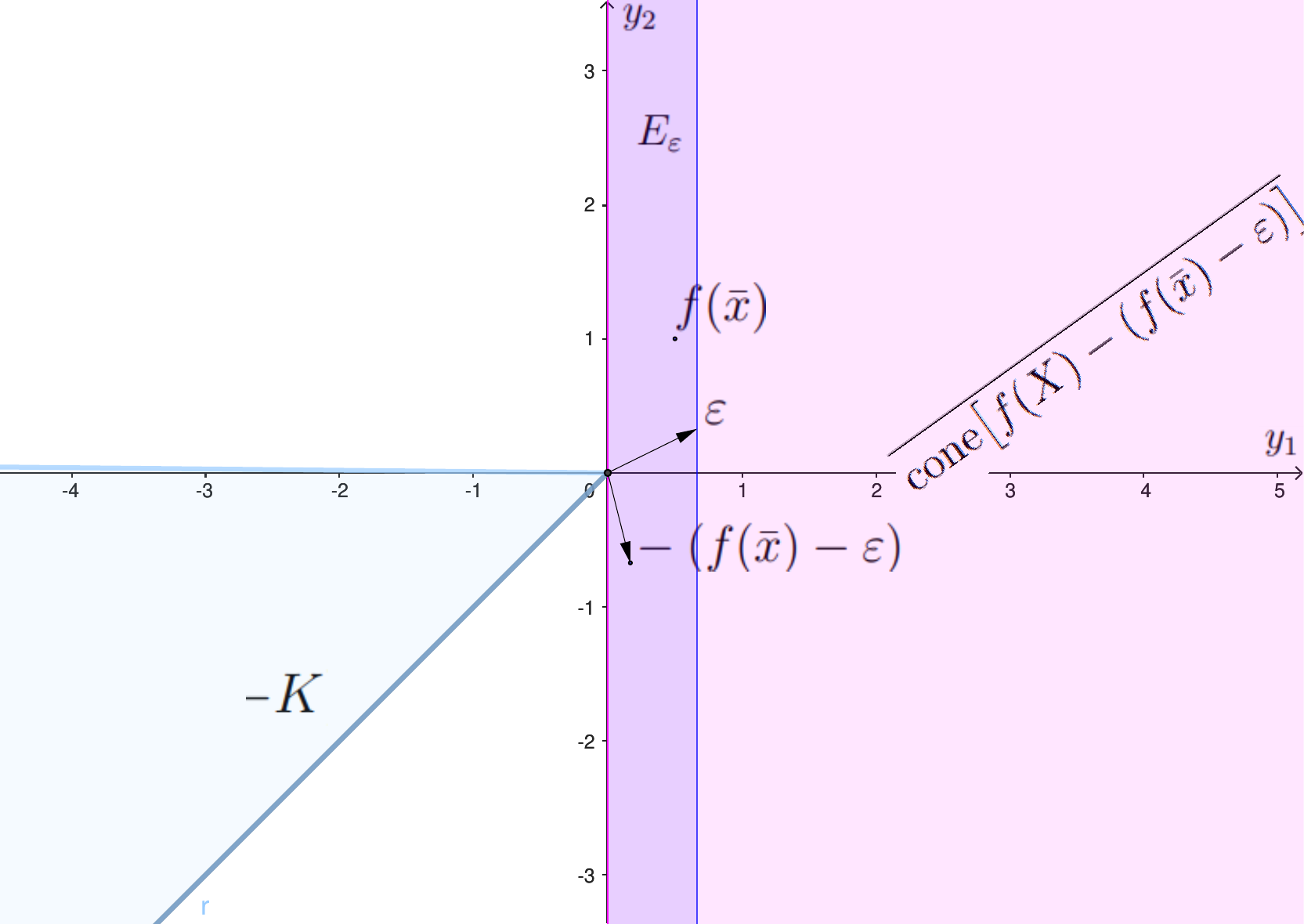}
				\caption{Illustration of Example \ref{T_Exam3}.}\label{Fig. 3}
			\end{center}
	\end{figure}} 
\end{Example}

\section{Approximate Proper Efficiency in Linear Vector Optimization}~\label{Sect_5}
\markboth{Approximate Proper Efficiency in Vector Optimization}{\sc n. t. t. huong}
\setcounter{equation}{0}

First, let us present the following classical result on proper efficiency in linear vector optimization. The detailed arguments in the proof can be effectively used for establishing our subsequent result.

\begin{Theorem}\label{thm_LVOP} {\rm (See, e.g.,~\cite[Theorem~3.4.7]{SNT_1985}}
	If {\rm (VP)} is a linear vector optimization problem with the polyhedral convex cone $K$ being pointed, then $E^{Be}=E$.
\end{Theorem}
\proof\  It suffices to show that $E\subset E^{Be}$. Let $\bar x\in E$ be given arbitrarily. Since the cone $K$ is pointed, by Remark~\ref{Remark} one has 
\begin{eqnarray}\label{eff_LVOP}
	{\rm cone}[f(X)+K-f(\bar x)]\cap (-K)=\{0\}.
\end{eqnarray}
Since the function $f:\mathbb R^n\to\mathbb R^m$ is an affine operator and $X$ is a polyhedral convex set, by~\cite[Theorems~19.1 and~19.3]{Roc70} we can assert that $f(X)$ is a polyhedral convex set. Then, since $K$ is a polyhedral convex cone, we can apply~\cite[Corollary~19.3.2]{Roc70} to obtain that $f(X)+K$ is a polyhedral convex set. Then, it follows that $f(X)+K-f(\bar x)$ is a polyhedral convex set. Therefore, as $0\in f(X)+K-f(\bar x)$, Corollary~19.7.1 from~\cite{Roc70} assures that ${\rm cone}\big[f(X)+K-f(\bar x)\big]$ is a polyhedral convex cone. In particular, this cone is closed. Consequently, the condition~\eqref{eff_LVOP} yields
\begin{equation*}\label{Be_LVOP} 
	\overline{\rm cone}\big[f(X)+K-f(\bar x)\big]\cap \big({-K}\big)=\{0\},
\end{equation*}
which justifies the inclusion $\bar x\in E^{Be}$.

The proof is complete. $\hfill\Box$

The two following statements will be used for the proof of Theorem~\ref{thm_LVOP_varepsilon_e}, which is the main result of this section.
\begin{Lemma}\label{recession} {\rm (See~\cite[Corollary 9.1.2]{Roc70})} 
	Let $D_1$ and $D_2$ be nonempty closed convex sets in $\mathbb R^m$. Assume there is no direction of recession of $D_1$ whose opposite is a direction of recession of $D_2$. (This is true in particular if either $D_1$ or $D_2$ is bounded.) Then $D_1+D_2$ is closed, and \, $0^+(D_1+D_2)=0^+D_1+0^+D_2.$
\end{Lemma}

\begin{Lemma}\label{rep poly conv cone} {\rm (See~\cite[Theorem 19.7]{Roc70})} 
	Let $C$ be a nonempty polyhedral convex set in $\mathbb R^m$, and let $D$ be the closure of the convex cone generated by $C$.
	Then $D$ is a polyhedral convex cone, and
	$$D=\left\{w \,:\, w=\lambda y\;, \lambda>0,\; y\in C \ \; {\rm or} \ \; w\in 0^+C\right\}.$$
\end{Lemma}

\begin{Theorem}\label{thm_LVOP_varepsilon_e} Let $e\in K\setminus\{0\}$. If {\rm (VP)} is a linear vector optimization problem with the polyhedral convex cone $K$ being pointed, then either $E^{Be}_e$ is empty, or $E^{Be}_e=E_e$.
\end{Theorem}
\proof\  Suppose that there exists $\bar x\in E^{Be}_e$. Since $K$ is pointed, by Definition~\ref{Benson'appro pro_e} and Lemma~\ref{equivalence}, the last inclusion means that
\begin{equation}\label{rep1} 
	\overline{\rm cone}\big[f(X)-(f(\bar x)-e)\big]\cap \big({-K}\big)=\{0\}.
\end{equation}
To prove the equality $E^{Be}_e=E_e$, we need to show that the equality 
\begin{equation}\label{rep2} 
	\overline{\rm cone}\big[f(X)-(f(\hat x)-e)\big]\cap \big({-K}\big)=\{0\}
\end{equation} holds for every $\hat x\in E_e$. 
Arguing by contradiction, suppose that there exists $\hat x\in E_e$ but~\eqref{rep2} fails to hold. Then we would find a nonzero vector $w\in -K$ with $$w\in \overline{\rm cone}[f(X)-(f(\hat x)-e)].$$ Since $f:\mathbb R^n\to\mathbb R^m$ is an affine operator and $X$ is a polyhedral convex set, by~\cite[Theorems~19.1 and~19.3]{Roc70} we can assert that $f(X)-(f(\hat x)-e)$ is a polyhedral convex set. Thus, by Lemma~\ref{rep poly conv cone}, either $w=\lambda[f(x)-(f(\hat x)-e)]$ for some $x\in X$ and $\lambda>0$ or  $w\in 0^+[f(X)-(f(\hat x)-e)]$.

First, suppose that $w=\lambda[f(x)-(f(\hat x)-e)]$ where $x\in X$ and $\lambda>0$. As $w\in{-K}$ and $w\ne 0$, one has $\lambda[f(x)-(f(\hat x)-e)]\in{-K}$ and $\lambda [f(x)-(f(\hat x)-e)]\ne 0$. This yields $\lambda[(f(\hat x)-e)-f(x)]\in K$ and $\lambda [f(x)-(f(\hat x)-e)]\ne 0$. Since $K$ is a cone, it follows that $(f(\hat x)-e)-f(x)\in K$ and $f(x)-(f(\hat x)-e)\ne 0$. Therefore, $$f(x)\leq_K f(\hat x)-e\ \; {\rm and}\ \; f(x)\neq f(\hat x)-e.$$ This contradicts the
assumption $\hat x\in E_e$.

Now, suppose that $w\in 0^+[f(X)-(f(\hat x)-e)]$. As the set $\{-(f(\hat x)-e)\}$ is bounded, by Lemma~\ref{recession} we have
$$\begin{array}{rl} 0^+[f(X)-(f(\hat x)-\varepsilon)]&=0^+[f(X)]+0^+[\{-(f(\hat x)-e)\}]\\ &=0^+[f(X)]+\{0\}\\&=0^+[f(X)].\end{array}$$
It follows that $$w\in 0^+[f(X)]= 0^+[f(X)-(f(\bar x)-e)]\subset \overline{\rm cone}\big[f(X)-(f(\bar x)-e)\big],$$ where the last inclusion is assured by Lemma~\ref{rep poly conv cone} and the convex polyhedrality of the set $f(X)-(f(\bar x)-e)$. The inclusion $w\in  \overline{\rm cone}\big[f(X)-(f(\bar x)-e)\big]$ contradicts~\eqref{rep1}. 

We have thus proved that~\eqref{rep2} holds for all $\hat x\in E_e$. The proof is complete. $\hfill\Box$

\medskip
Applying Theorem \ref{thm_LVOP_varepsilon_e} for $K=\mathbb{R}^m_+$, we have the next result.

\begin{Theorem}\label{LVOP_varepsilon}  Let $\varepsilon\in\mathbb{R}^m_+\setminus\{0\}$. If {\rm (VP)} is a linear vector optimization problem with the ordering cone $\mathbb{R}^m_+$, then either $E^{Be}_\varepsilon$ is empty, or $E^{Be}_\varepsilon=E_\varepsilon$.
\end{Theorem}

\noindent
\textbf{Acknowledgements}  The author would to thank Professor Nguyen Dong Yen for helpful discussions on the subject.


\begin{thebibliography}{99}%Please arrange author's name in alphabetical order.
	
\bibitem{Be79} B. Benson, \textit{An improved definition of proper efficiency for vector maximization with respect to cones}, J. Math. Anal. Appl. 71 (1979), 232--241. 
	
\bibitem{Bo77}  J. M. Borwein, \textit{Proper efficient points for maximizations with respect to cones}, SIAM J. Control Optim.~15 (1977), 57--63 
	
\bibitem{Choo84} E. U. Choo, \textit{Proper efficiency and the linear fractional vector maximum problem}, Oper. Res.~32 (1984), 216--220.

\bibitem{CK2016} T. D. Chuong, D. S. Kim, \textit{Approximate solutions of multiobjective optimization problems}, Positivity~20 (2016), 187–207. 

\bibitem{GE68} A. M. Geoffrion, \textit{Proper efficiency and the theory of vector maximization}, J. Math. Anal. Appl.~22 (1968), 613--630. 

\bibitem{GMZ94} A. Guerraggio, E. Molho, A. Zaffaroni,  \textit{On the notion of proper efficiency in vector optimization}. J. Optim. Theory Appl. 82 (1994), 1--21.

\bibitem{GJN06} C. Guti\'errez, B. Jimenez, V. Novo, \textit{A unified approach and optimality conditions for approximate solutions of vector optimization problems}, SIAM J. Optim.~17 (2006), 688--710.

\bibitem{Gutierrez-et-al-12} C. Guti\'errez, L. Huerga, V. Novo, \textit{Scalarization and saddle points of approximate proper
solutions in nearly subconvexlike vector optimization problems}, J. Math. Anal. Appl.~389 (2012), 1046–1058.

\bibitem{Henig_JOTA1982}  M. I. Henig, \textit{Proper efficiency with respect to cones}, J. Optim. Theory Appl.~36 (1982), 387--407.
	
\bibitem{HWYY2019} N. T. T. Huong, C.-F. Wen,  J.-C. Yao,  N. D. Yen, \textit{Proper efficiency in linear fractional vector optimization via Benson's characterization}, Optimization~72 (2023), 263–276. 

\bibitem{HYY_JOGO2020} N. T. T. Huong,  J.-C. Yao, N. D. Yen, \textit{Geoffrion's proper efficiency in linear fractional vector optimization with unbounded constraint sets}, J. Global Optim.~78 (2020), 545--562.
	
\bibitem{HYY2020} N. T. T. Huong, J.-C. Yao,  N. D. Yen, \textit{New results on proper efficiency for a class of vector optimization problems}, Appl. Anal.~15 (2021), 3199–3211.

\bibitem{HY2022a} N. T. T. Huong, N. D. Yen, \textit{Improperly efficient solutions in a class of vector optimization problems}, J. Global Optim.~82 (2022), 375--387. 

\bibitem{I74} H. Isermann \textit{ Proper efficiency and the linear vector maximum problem}, Operations Res.~22 (1974), 189--191.

%\bibitem{KMST2021} D. S. Kim, B. S. Mordukhovich, T. S. Pham, N. V. Tuyen \textit{Existence of efficient and properly efficient solutions to problems of constrained vector optimization}, Math. Program~190 (2021), 259--283.

\bibitem{LTY05}  G. M. Lee, N. N. Tam,  N. D. Yen, \textit{Quadratic Programming and Affine Variational Inequalities: A Qualitative Study}, Springer Verlag, New York, 2005.
	
\bibitem{LiWang} Z. Li, S. Wang, \textit {$\varepsilon$-approximate solutions in multiobjective optimization}, Optimization~44 (1998), 161--174.
	
\bibitem{Liu99}  J.-C. Liu, \textit {$\varepsilon$-properly efficient solutions to nondifferentiable multiobjective programming problems}, Appl. Math. Lett.~12 (1999), 109--113.

\bibitem{Luc89}  D. T. Luc, \textit {Theory of Vector Optimization} Springer, Berlin (1989).

%\bibitem{Nieuwenhuis81} J. W. Nieuwenhuis, \textit{Properly efficient and efficient solutions for vector maximization problems
%in Euclidean space}, J. Math. Anal. Appl., 84 (1981), pp. 311--317.
	
\bibitem{Roc70} R. T. Rockafellar, \textit{Convex Analysis}, Princeton University Press, Princeton, New Jersey, 1970. 

\bibitem{RM2000} W. D. Rong, Y. Ma, \textit{e-properly efficient solutions of vector optimization problems with set-valued
maps}, OR Trans.~4 (2000), pp. 21--32.

\bibitem{SNT_1985} Y. Sawaragi, H. Nakayama, T. Tanino, \textit{Theory of Multiobjective Optimization}, Academic Press, Inc., Orlando, FL, 1985.
	
\bibitem{S86} R. E. Steuer, \textit{Multiple Criteria Optimization: Theory, computation and application}, John Wiley \& Sons, 1986, New York.

\bibitem{ST14} B. Soleimani, C. Tammer, \textit{Concepts for approximate solutions of vector optimization problems with variable order structures} Vietnam J. Math.~42 (2014), 543--566. 

\bibitem{T22} N. V. Tuyen, \textit{A note on approximate proper efficiency in linear fractional vector optimization}, Optim. Lett.~16 (2022), 1835–1845. 

\bibitem{Tammer94} C. Tammer, \textit{Stability results for approximately efficient solutions}, OR Spectrum~16 (1994),
47--52.

\bibitem{White86} D. J. White, \textit{Epsilon efficiency}, J. Optim. Theory Appl.~49 (1986), 319--337.

\bibitem{ZhaoYang15} K. Q. Zhao, X. M. Yang, \textit{E-Benson proper efficiency in vector optimization},  Optimization~64 (2015), 739--752.

\bibitem{ZhaoChenYang15} K. Q. Zhao, G. Y. Chen, X. M. Yang, \textit{Approximate proper efficiency in vector optimization}, Optimization~64 (2015), 1777--1793.

\end{thebibliography}
\end{document}